\documentstyle{article}
\begin{document}

\title{Bicartesian Coherence}
\author{K{\footnotesize OSTA} D{\footnotesize O\v SEN} and
Z{\footnotesize ORAN} P{\footnotesize ETRI\' C}\\[0.5cm]
Matemati\v cki institut \\
SANU, Knez Mihailova 35, p.f. 367 \\
11001 Belgrade, Serbia \\
email: \{kosta, zpetric\}@mi.sanu.ac.yu}
\date{}
\maketitle

\begin{abstract}
\noindent Coherence is demonstrated for categories with binary
products and sums, but without the terminal and the initial
object, and without distribution. This coherence amounts to the
existence of a faithful functor from a free category with binary
products and sums to the category of relations on finite ordinals.
This result is obtained with the help of proof-theoretic
normalizing techniques. When the terminal object is present,
coherence may still be proved if of binary sums we keep just their
bifunctorial properties. It is found that with the simplest
understanding of coherence this is the best one can hope for in
bicartesian categories. The coherence for categories with binary
products and sums provides an easy decision procedure for equality
of arrows. It is also used to demonstrate that the categories in
question are maximal, in the sense that in any such category that
is not a preorder all the equations between arrows involving only
binary products and sums are the same. This shows that the usual
notion of equivalence of proofs in nondistributive
conjunctive-disjunctive logic is optimally defined: further
assumptions would make this notion collapse into triviality. (A
proof of coherence for categories with binary products and sums
simpler than that presented in this paper may be found in Section
9.4 of {\it Proof-Theoretical Coherence}, revised version of
September 2007,
http://\-www.\-mi.\-sanu.\-ac.\-yu/\-$\sim$kosta/\-coh.\-pdf.)
\end{abstract}

\vspace{1cm}


  \vspace{1cm}

\section{Introduction}

At the very beginning of his categorial proof-theoretical program, Lambek
has formulated the idea that two proofs with the same premise and conclusion
should be considered equivalent iff they have the same ``generality'' (see
\cite{lam68}, p. 316, and \cite{lam69}, p. 89). The standard example of two
proofs with a different generality is given by the first projection and the second
projection from $p\wedge p$ to $p$; they generalize respectively to proofs
from $p\wedge q$ to $p$ and from $p\wedge q$ to $q$. Lambek's way of making
the idea of generality precise ran into difficulties (see \cite{lam72}, p.
65). Szabo pursued the idea in his own way in \cite{sza74} and \cite{sza89}.

We find that the most simple way to understand generality is to connect by a
link propositional letters that must remain identical after generalizing,
and not connect those that may differ. So for the first and the second
projection proof from $p\wedge p$ to $p$ we would have the two diagrams

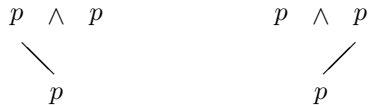
\begin{center}
\begin{picture}(150,70)
\put(25,50){\makebox(0,0){$\wedge$}}
\put(125,50){\makebox(0,0){$\wedge$}}
\put(10,50){\makebox(0,0){$p$}}
\put(40,50){\makebox(0,0){$p$}}
\put(110,50){\makebox(0,0){$p$}}
\put(140,50){\makebox(0,0){$p$}}
\put(25,20){\makebox(0,0){$p$}}
\put(125,20){\makebox(0,0){$p$}}
\put(12,40){\line(1,-1){12}}
\put(138,40){\line(-1,-1){12}}
\end{picture}
\end{center}

\noindent The two proofs would not be equivalent because these diagrams are different.

Such diagrams compose in an obvious way by composing links, and we have also
obvious identity diagrams like, for example,

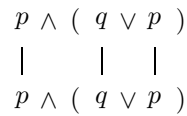
\begin{center}
\begin{picture}(80,70)
\put(20,50){\makebox(0,0){$\wedge$}}
\put(20,20){\makebox(0,0){$\wedge$}}
\put(50,50){\makebox(0,0){$\vee$}}
\put(50,20){\makebox(0,0){$\vee$}}
\put(30,50){\makebox(0,0){$($}}
\put(30,20){\makebox(0,0){$($}}
\put(70,50){\makebox(0,0){$)$}}
\put(70,20){\makebox(0,0){$)$}}
\put(10,50){\makebox(0,0){$p$}}
\put(40,50){\makebox(0,0){$q$}}
\put(60,50){\makebox(0,0){$p$}}
\put(10,20){\makebox(0,0){$p$}}
\put(40,20){\makebox(0,0){$q$}}
\put(60,20){\makebox(0,0){$p$}}
\put(10,40){\line(0,-1){10}}
\put(40,40){\line(0,-1){10}}
\put(60,40){\line(0,-1){10}}
\end{picture}
\end{center}

\noindent So diagrams, or a formal analogue of them, make a category. We call such
categories {\em graphical categories}. In this paper, the graphical category
will be the category of relations on finite ordinals.

Equivalence between proofs in intuitionistic logic is axiomatized
independently of these diagrams in the typed lambda calculus and in various
sorts of categories, like bicartesian closed categories. There, proofs are
coded by typed lambda terms or by arrow terms, and two proofs are considered
equivalent iff the coding terms are equal as lambda terms or as arrow terms
in categories. This approach is rather standard nowadays, because lambda
equality and equality of arrows in categories match pretty well the
equivalence between proofs induced by normalization in natural deduction or
cut elimination in sequent systems for intuitionistic logic.

The question then arizes how this standard notion of equivalence relates to
generality of proofs. This question can be exactly posed by asking whether
for a freely generated category ${\cal C}$ of some sort, like, for example,
a freely generated bicartesian closed category, there is a faithful functor $%
G$ from ${\cal C}$ to a graphical category ${\cal G}$, in which diagrams
like the two diagrams above, or a formal analogue of these diagrams, would
be arrows. In limit cases, the functor $G$ may be an isomorphism, but in
general it is enough that it be faithful, which means that for two arrow
terms $f,g:A\rightarrow B$ of ${\cal C}$ we shall have

\begin{center}
$f=g$ {\em in} ${\cal C}$ \quad {\em iff} \quad $G(f)=G(g)$ {\em in} ${\cal G%
}$.
\end{center}

\noindent From left to right, this equivalence follows from $G$'s being a
functor, and from right to left it expresses the faithfulness of $G$.

Understood in this manner, Lambek's generality idea merges with coherence
questions in category theory, which were treated intensively when Lambek
wrote the three papers cited above. Now, coherence in category theory has no
doubt been understood in various ways. (We shall not try to survey here the
literature on this question; for earlier works see \cite{mac72}.) Although
Mac Lane's paradoxical dictum ``All diagrams commute'' can be made precise
in different ways, the paradigmatic results on coherence of \cite{mac63} and
\cite{kel71} can be understood as faithfulness results, like the equivalence
above, and this is how we understand coherence here. The graphical
categories of these results, whose investigation starts with \cite{eil66},
are very well adapted to make precise Lambek's idea of generality in logic.

The faithfulness equivalence above, which may also be understood as a
coherence equivalence, is, from a logical point of view, a completeness
equivalence. The freely generated category ${\cal C}$ is syntax, i.e. a
formal system, the graphical category ${\cal G}$ is a model, the
functoriality of $G$, i.e. the implication from left to right, is soundness
(here it is desirable that $G$ preserve also the particular structure of $%
{\cal C}$, and not only identities and composition), and the implication
from right to left is completeness proper. For this completeness result to
be interesting, there should be a gain in passing from ${\cal C}$ to ${\cal G%
}$. It is desirable that ${\cal G}$ be easy to handle, so that, for example,
we may decide equality of arrows in ${\cal C}$ by passing to ${\cal G}$, or
that we may normalize arrow terms by referring to ${\cal G}$, without going
through tedious syntactic reductions.

If we understand generality in ${\cal G}$ in the simplest fashion: ``Connect
by a link all propositional letters that must remain the same after
generalizing'', the generality idea matches equality in categories that
correspond to intuitionistic propositional logic only to a limited extent.
(It fares better in linear logic, as the coherence result of \cite{kel71}
shows.) We have coherence, i.e. faithfulness, i.e. soundness and
completeness, for categories with binary product, which covers the purely
conjunctive fragment of logic, and for cartesian categories, where the
terminal object is added, which covers the conjunctive fragment extended
with the constant true proposition (see references in Section 4). By
duality, this covers also the purely disjunctive fragment, and this fragment
extended with the constant absurd proposition.

We shall see in this paper that we have coherence also for categories with
binary products and sums (i.e. coproducts), which covers the
conjunctive-disjunctive fragment without distribution of conjunction over
disjunction. We shall also see that we have coherence for cartesian
categories, where the terminal object is present, extended with an operation
on objects and arrows that keeps of sum (or of product) just its
bifunctorial properties.

Coherence fails for bicartesian categories, where the terminal and the
initial object have been added. If $\top$ is the terminal object, which
corresponds to the constant true proposition, and $\bot$ the initial object,
which corresponds to the constant absurd proposition, both for the first and
for the second projection proof from $\bot \wedge \bot$ to $\bot$ we have
the empty diagram, and analogously for the first and the second injection
proof from $\top$ to $\top\vee\top$. However, neither the first and the
second projection from $\bot \wedge \bot$ to $\bot$, nor the first and the
second injection from $\top$ to $\top\vee\top$, are equal in all bicartesian
categories (see Section 4). While the first and the second projection from $%
\bot\wedge\bot$ to $\bot$ become equal in all bicartesian closed categories,
the first and the second injection from $\top$ to $\top\vee\top$ are equal
in a bicartesian closed category iff the category has collapsed into a
preorder, i.e. a category where all arrows with the same source and the same
target are equal. (The equality between the first and the second projection
from $\bot\wedge\bot$ to $\bot$ in bicartesian closed categories, which is a
consequence of the existence of a right adjoint to the functor $\bot \times$%
, implies that in bicartesian closed categories all arrows with the same
source and the target $\bot$ are equal; cf. Section 5 and \cite{lam86}, p.
67, Proposition 8.3. The equality between the first and the second injection
from $\top$ to $\top\vee\top$ yields preordering in bicartesian closed
categories because in these categories all arrows are in one-to-one
correspondence with arrows whose source is $\top$; cf. Section 5.) It is the
completeness part of coherence that fails for bicartesian categories---the
soundness part holds true.

It is rather typical for special objects to make trouble for coherence
results. Such was the case for the unit object in symmetric monoidal closed
categories (see \cite{kel71}), and such is the case in bicartesian
categories with the terminal and the initial object. (Special objects may
also cause trouble for normalizing terms; see, for example, \cite{lam86}, p.
88.) So, for some purposes, it may be unwise to subsume the terminal and the
initial object under a generalized concept of finite product and sum, where
they are the nullary cases. This may obscure matters.

Even the soundness part of coherence fails for distributive categories with
binary products and sums, and for distributive bicartesian categories (for
these categories see \cite{law97}, pp. 222-223 and Session 26, and \cite
{coc93}) . If we have a distributivity isomorphism from $p\wedge (q \vee r)$
to $(p\wedge q) \vee (p\wedge r)$, then the composition of diagrams

\begin{center}
\begin{picture}(120,80)
\put(10,70){\makebox(0,0){$($}}
\put(20,70){\makebox(0,0){$p$}}
\put(30,70){\makebox(0,0){$\wedge$}}
\put(40,70){\makebox(0,0){$q$}}
\put(50,70){\makebox(0,0){$)$}}
\put(60,70){\makebox(0,0){$\vee$}}
\put(70,70){\makebox(0,0){$($}}
\put(80,70){\makebox(0,0){$p$}}
\put(90,70){\makebox(0,0){$\wedge$}}
\put(100,70){\makebox(0,0){$r$}}
\put(110,70){\makebox(0,0){$)$}}

\put(40,40){\makebox(0,0){$p$}}
\put(50,40){\makebox(0,0){$\wedge$}}
\put(60,40){\makebox(0,0){$($}}
\put(70,40){\makebox(0,0){$q$}}
\put(80,40){\makebox(0,0){$\vee$}}
\put(90,40){\makebox(0,0){$r$}}
\put(100,40){\makebox(0,0){$)$}}

\put(10,10){\makebox(0,0){$($}}
\put(20,10){\makebox(0,0){$p$}}
\put(30,10){\makebox(0,0){$\wedge$}}
\put(40,10){\makebox(0,0){$q$}}
\put(50,10){\makebox(0,0){$)$}}
\put(60,10){\makebox(0,0){$\vee$}}
\put(70,10){\makebox(0,0){$($}}
\put(80,10){\makebox(0,0){$p$}}
\put(90,10){\makebox(0,0){$\wedge$}}
\put(100,10){\makebox(0,0){$r$}}
\put(110,10){\makebox(0,0){$)$}}

\put(20,65){\line(1,-1){18}}
\put(40,65){\line(3,-2){28}}
\put(80,65){\line(-2,-1){38}}
\put(100,65){\line(-1,-2){10}}
\put(38,35){\line(-1,-1){18}}
\put(42,35){\line(2,-1){38}}
\put(70,35){\line(-3,-2){28}}
\put(90,35){\line(1,-2){10}}
\end{picture}
\end{center}

\noindent  yields the diagram on the left-hand side below,
which doesn't amount to the
identity diagram on the right-hand side:

\begin{center}
\begin{picture}(320,50)
\put(10,40){\makebox(0,0){$($}}
\put(20,40){\makebox(0,0){$p$}}
\put(30,40){\makebox(0,0){$\wedge$}}
\put(40,40){\makebox(0,0){$q$}}
\put(50,40){\makebox(0,0){$)$}}
\put(60,40){\makebox(0,0){$\vee$}}
\put(70,40){\makebox(0,0){$($}}
\put(80,40){\makebox(0,0){$p$}}
\put(90,40){\makebox(0,0){$\wedge$}}
\put(100,40){\makebox(0,0){$r$}}
\put(110,40){\makebox(0,0){$)$}}

\put(10,10){\makebox(0,0){$($}}
\put(20,10){\makebox(0,0){$p$}}
\put(30,10){\makebox(0,0){$\wedge$}}
\put(40,10){\makebox(0,0){$q$}}
\put(50,10){\makebox(0,0){$)$}}
\put(60,10){\makebox(0,0){$\vee$}}
\put(70,10){\makebox(0,0){$($}}
\put(80,10){\makebox(0,0){$p$}}
\put(90,10){\makebox(0,0){$\wedge$}}
\put(100,10){\makebox(0,0){$r$}}
\put(110,10){\makebox(0,0){$)$}}

\put(18,35){\line(0,-1){20}}
\put(40,35){\line(0,-1){20}}
\put(82,35){\line(0,-1){20}}
\put(100,35){\line(0,-1){20}}
\put(22,35){\line(3,-1){56}}
\put(78,35){\line(-3,-1){56}}

\put(210,40){\makebox(0,0){$($}}
\put(220,40){\makebox(0,0){$p$}}
\put(230,40){\makebox(0,0){$\wedge$}}
\put(240,40){\makebox(0,0){$q$}}
\put(250,40){\makebox(0,0){$)$}}
\put(260,40){\makebox(0,0){$\vee$}}
\put(270,40){\makebox(0,0){$($}}
\put(280,40){\makebox(0,0){$p$}}
\put(290,40){\makebox(0,0){$\wedge$}}
\put(300,40){\makebox(0,0){$r$}}
\put(310,40){\makebox(0,0){$)$}}
\put(210,10){\makebox(0,0){$($}}
\put(220,10){\makebox(0,0){$p$}}
\put(230,10){\makebox(0,0){$\wedge$}}
\put(240,10){\makebox(0,0){$q$}}
\put(250,10){\makebox(0,0){$)$}}
\put(260,10){\makebox(0,0){$\vee$}}
\put(270,10){\makebox(0,0){$($}}
\put(280,10){\makebox(0,0){$p$}}
\put(290,10){\makebox(0,0){$\wedge$}}
\put(300,10){\makebox(0,0){$r$}}
\put(310,10){\makebox(0,0){$)$}}

\put(220,35){\line(0,-1){20}}
\put(240,35){\line(0,-1){20}}
\put(280,35){\line(0,-1){20}}
\put(300,35){\line(0,-1){20}}
\end{picture}
\end{center}

In distributive bicartesian categories we also have that $p \wedge \bot$ is
isomorphic to $\bot$, but the composition of two empty diagrams for this
isomorphism is not equal to the identity diagram from $p\wedge \bot$ to $%
p\wedge \bot$.

Without disjunction, but with implication, the situation is not better. Both
the soundness part and the completeness part of coherence fail for cartesian
closed categories. For soundness, we have the counterexample

\begin{center}
\begin{picture}(400,80)(25,0)
\put(10,10){\makebox(0,0){$($}}
\put(20,10){\makebox(0,0){$q$}}
\put(30,10){\makebox(0,0){$\rightarrow$}}
\put(40,10){\makebox(0,0){$($}}
\put(50,10){\makebox(0,0){$q$}}
\put(60,10){\makebox(0,0){$\wedge$}}
\put(70,10){\makebox(0,0){$p$}}
\put(80,10){\makebox(0,0){$)$}}
\put(90,10){\makebox(0,0){$)$}}
\put(100,10){\makebox(0,0){$\wedge$}}
\put(110,10){\makebox(0,0){$($}}
\put(120,10){\makebox(0,0){$q$}}
\put(130,10){\makebox(0,0){$\rightarrow$}}
\put(140,10){\makebox(0,0){$($}}
\put(150,10){\makebox(0,0){$q$}}
\put(160,10){\makebox(0,0){$\wedge$}}
\put(170,10){\makebox(0,0){$p$}}
\put(180,10){\makebox(0,0){$)$}}
\put(190,10){\makebox(0,0){$)$}}

\put(60,40){\makebox(0,0){$q$}}
\put(70,40){\makebox(0,0){$\rightarrow$}}
\put(80,40){\makebox(0,0){$($}}
\put(90,40){\makebox(0,0){$q$}}
\put(100,40){\makebox(0,0){$\wedge$}}
\put(110,40){\makebox(0,0){$p$}}
\put(120,40){\makebox(0,0){$)$}}

\put(100,70){\makebox(0,0){$p$}}

\put(100,65){\line(1,-2){10}}
\put(75,45){\oval(30,16)[t]}

\put(58,35){\line(-2,-1){38}}
\put(88,35){\line(-2,-1){38}}
\put(108,35){\line(-2,-1){38}}
\put(62,35){\line(3,-1){58}}
\put(92,35){\line(3,-1){58}}
\put(112,35){\line(3,-1){58}}

\put(210,10){\makebox(0,0){$($}}
\put(220,10){\makebox(0,0){$q$}}
\put(230,10){\makebox(0,0){$\rightarrow$}}
\put(240,10){\makebox(0,0){$($}}
\put(250,10){\makebox(0,0){$q$}}
\put(260,10){\makebox(0,0){$\wedge$}}
\put(270,10){\makebox(0,0){$p$}}
\put(280,10){\makebox(0,0){$)$}}
\put(290,10){\makebox(0,0){$)$}}
\put(300,10){\makebox(0,0){$\wedge$}}
\put(310,10){\makebox(0,0){$($}}
\put(320,10){\makebox(0,0){$q$}}
\put(330,10){\makebox(0,0){$\rightarrow$}}
\put(340,10){\makebox(0,0){$($}}
\put(350,10){\makebox(0,0){$q$}}
\put(360,10){\makebox(0,0){$\wedge$}}
\put(370,10){\makebox(0,0){$p$}}
\put(380,10){\makebox(0,0){$)$}}
\put(390,10){\makebox(0,0){$)$}}

\put(290,40){\makebox(0,0){$p$}}
\put(310,40){\makebox(0,0){$\wedge$}}
\put(330,40){\makebox(0,0){$p$}}

\put(310,70){\makebox(0,0){$p$}}

\put(308,65){\line(-1,-1){18}}
\put(312,65){\line(1,-1){18}}

\put(235,15){\oval(30,16)[t]}
\put(335,15){\oval(30,16)[t]}

\put(290,35){\line(-1,-1){20}}
\put(330,35){\line(2,-1){40}}

\end{picture}
\end{center}

\noindent while for completeness, counterexamples may be constructed along
the lines
of \cite{sza75}. (Both soundness and completeness would fail if we had only
implication, but the associated categories, which correspond to the lambda
calculus without product types, are not usually considered.)

So it seems that with the results of this paper we have reached the limits
of coherence for the simple approach to generality in logic. This doesn't
mean that another approach, with a more subtly built graphical category $%
{\cal G}$, couldn't vindicate Lambek's idea in wider fragments of logic.

The main coherence result we are going to establish here is useful for a
particular purpose. Recently, Cockett and Seely set forth in \cite{coc00} a
syntactical decision procedure obtained via cut elimination for equality of
arrow terms in freely generated categories with finite products and sums,
without distribution, which amount to bicartesian categories\footnote{%
We are grateful to Robert Seely and Robin Cockett for having sent us the
draft of \cite{coc00}, which prompted the writing of this paper.}. (That
there is such a procedure is claimed in \cite{sza78}, p. 69, but we have
found that difficult to check.) We don't cover the whole ground of \cite
{coc00}, since coherence fails iff we add the empty product and the empty
sum, i.e. the terminal and the initial object. However, as far as it goes,
for equations between terms involving only nonempty finite products and
sums, which covers the major part of the categories in question, our
coherence yields a very simple graphical decision procedure, in which we
find a clear advantage over syntactical procedures, even when they are
entirely explicit.

We use our coherence result for categories with binary products and sums to
demonstrate in the last section of this paper that these categories are
maximal, in the sense that if such a category satisfies every instance of
any equation not satisfied in freely generated categories of this sort, then
this category is a preorder. Analogous results hold for cartesian and
cartesian closed categories (see \cite{dos00}, \cite{sim95} and \cite{dos00a}%
). Such results are interesting for logic, because they show that our choice
of equations is optimal. These equations are wanted, because they are
induced by normalization of proofs, and no equation is missing, because any
further equation would lead to collapse: all arrows with the same source and
the same target would be equal. Our maximality result in the last section
shows optimal the choice of equations assumed for conjunctive-disjunctive
logic in the absence of distribution.

The literature on bicartesian categories, without distribution, and without
closure, i.e. exponentiation, does not seem very rich, though, of course,
the notions of product and sum (coproduct) are explained in every textbook
of category theory. An early reference we know about, but which does not
cover matters we are treating, is \cite{eck62}. Of course, nondistributive
lattices have been extensively studied, but though this topic is related to
bicartesian categories, categorial studies are on a different level. In the
literature on nonclassical logics, one encounters very much studied logics
where conjunction and disjunction make a nondistributive lattice
structure---such is, for instance, linear logic---, but the purely
conjunctive-disjunctive fragment is not usually separated and considered
categorially. (A sequent system for nondistributive conjunction and
disjunction is considered towards the end of \cite{reb93} from the point of
view of algebraic logic, but categories are not mentioned there.)

\section{Free Bicartesian Categories}

The propositional language ${\cal P}$ is generated from a set of {\em %
propositional letters} ${\cal L}$ with the nullary connectives, i.e.
propositional constants, I and O, and the binary connectives $\times$ and $+$%
. The fragments ${\cal P}_{\times,+,\mbox{\scriptsize{\rm I}}}$, ${\cal P}%
_{\times,+}$ etc. of ${\cal P}$ are obtained by keeping only those formulae
of ${\cal P}$ that contain the connectives in the index. For the
propositional letters of ${\cal P}$, i.e. for the members of ${\cal L}$, we
use the schematic letters $p,q,\ldots,p_1,\ldots,$ and for the formulae of $%
{\cal P}$, or of its fragments, we use the schematic letters $%
A,B,\ldots,A_1,\ldots$

Next we define inductively the {\em terms} that will stand for the arrows of
the free bicartesian category ${\cal C}$ generated by ${\cal L}$. Every term
has a {\em type}, which is a pair $(A,B)$ of formulae of ${\cal P}$. That a
term $f$ is of type $(A,B)$ is written $f:A\rightarrow B$. The {\em atomic}
terms of ${\cal C}$ are for every $A$ of ${\cal P}$
\[
\begin{array}{ccc}
& \mbox{\bf 1}_A:A\rightarrow A, &  \\
k_A:A\rightarrow \mbox{\rm I}, &  & l_A:\mbox{\rm O}\rightarrow A.
\end{array}
\]
The terms $\mbox{\bf 1}_A$ are called {\em identities}. The other terms of $%
{\cal C}$ are generated with the following operations on terms, which we
present by rules so that from the terms in the premises we obtain the terms
in the conclusion:
\[
\frac{f:A\rightarrow B \quad g:B\rightarrow C}{g\circ f:A \rightarrow C}
\]
\[
\begin{array}{cc}
\displaystyle\frac{f:A\rightarrow C}{K^1_B f:A\times B\rightarrow C} & %
\displaystyle\frac{f:C\rightarrow A}{L^1_B f:C\rightarrow A+B} \\[.3cm]
\displaystyle\frac{f:B\rightarrow C}{K^2_A f:A\times B\rightarrow C} & %
\displaystyle\frac{f:C\rightarrow B}{L^2_A f:C\rightarrow A+B} \\[.3cm]
\displaystyle\frac{f:C\rightarrow A \quad g:C\rightarrow B}{\langle
f,g\rangle:C\rightarrow A\times B} & \displaystyle\frac{f:A\rightarrow C
\quad g:B\rightarrow C}{[ f,g ]:A+B\rightarrow C}
\end{array}
\]

\noindent We use $f,g,\ldots,f_1,\ldots$ as schematic letters for terms of $%
{\cal C}$.

The category ${\cal C}$ has as objects the formulae of ${\cal P}$ and as
arrows equivalence classes of terms so that the following equations are
satisfied for $i\in\{1,2\}$:
\[
\begin{array}{l}
(cat\: 1)\quad \mbox{\bf 1}_B\circ f=f\circ \mbox{\bf 1}_A=f, \\
(cat\: 2)\quad h\circ (g\circ f)=(h\circ g)\circ f,
\end{array}
\]
\[
\begin{array}{ll}
(K1)\quad g\circ K^i_A f=K^i_A(g\circ f), & (L1)\quad L^i_A g \circ
f=L^i_A(g\circ f), \\[.05cm]
(K2)\quad K^i_Ag\circ\langle f_1,f_2\rangle =g\circ f_i, & (L2)\quad [g_1,
g_2]\circ L^i_A f=g_i\circ f, \\[.05cm]
(K3)\quad \langle g_1,g_2\rangle \circ f=\langle g_1\circ f, g_2\circ
f\rangle, & (L3)\quad g\circ[f_1,f_2]=[g\circ f_1, g\circ f_2], \\[.05cm]
(K4)\quad \langle K^1_B \mbox{\bf 1}_A, K^2_A \mbox{\bf 1}_B \rangle=%
\mbox{\bf 1}_{A\times B}, & (L4)\quad [L^1_B \mbox{\bf 1}_A, L^2_A%
\mbox{\bf
1}_B]=\mbox{\bf 1}_{A+B}, \\[.05cm]
(k)\quad {\mbox{\rm for }} f:A\rightarrow \mbox{\rm I},\; f=k_A, & (l)\quad {%
\mbox{\rm for }} f:\mbox{\rm O}\rightarrow A,\; f=l_A.
\end{array}
\]

A category ${\cal C}^{\prime}$ isomorphic to ${\cal C}$ is obtained with the
same objects, and terms defined inductively as follows. The atomic terms are
for every $A$ and every $B$ of ${\cal P}$
\[
\begin{array}{lcl}
& \mbox{\bf 1}_A:A\rightarrow A, &  \\
k_A:A\rightarrow\mbox{\rm I}, &  & l_A:\mbox{\rm O}\rightarrow A, \\
k^1_{A,B}:A\times B\rightarrow A, &  & l^1_{A,B}:A\rightarrow A+B, \\
k^2_{A,B}:A\times B\rightarrow B, &  & l^2_{A,B}:B\rightarrow A+B, \\
w_A:A\rightarrow A\times A, &  & m_A:A+A\rightarrow A,
\end{array}
\]
and we have the following operations on terms:
\[
\frac{f:A\rightarrow B \quad g:B\rightarrow C}{g\circ f: A\rightarrow C}
\]
\[
\frac{f:A\rightarrow B \quad g:C\rightarrow D}{f\times g:A\times
C\rightarrow B\times D} \quad\quad\quad \frac{f:A\rightarrow B \quad
g:C\rightarrow D}{f+g:A+C\rightarrow B+D}
\]

On these terms we impose the equations $(cat\: 1)$, $(cat\: 2)$, $(k)$, $(l)$
and
\[
\begin{array}{l}
(\times 1)\; \mbox{\bf 1}_A\times\mbox{\bf 1}_B=\mbox{\bf 1}_{A\times B}, \\%
[.05cm]
(\times 2)\; (g_1\circ g_2)\times(f_1\circ f_2)= (g_1\times
f_1)\circ(g_2\times f_2), \\[.05cm]
(k^i)\; k^i_{B_1,B_2}\circ(f_1\times f_2)=f_i\circ k^i_{A_1,A_2}, \\[.05cm]
(w)\; w_B\circ f=(f\times f)\circ w_A, \\[.05cm]
(kw1)\; k^i_{A,A}\circ w_A=\mbox{\bf 1}_A, \\[.05cm]
(kw2)\; (k^1_{A,B}\times k^2_{A,B})\circ w_{A\times B}= \mbox{\bf 1}%
_{A\times B},
\end{array}
\]
\[
\begin{array}{l}
(+1)\; \mbox{\bf 1}_A+\mbox{\bf 1}_B=\mbox{\bf 1}_{A+B}, \\[.05cm]
(+2)\; (g_1\circ g_2)+(f_1\circ f_2)=(g_1+f_1)\circ(g_2+f_2), \\[.05cm]
(l^i)\; (f_1+f_2)\circ l^i_{A_1,A_2}=l^i_{B_1,B_2}\circ f_i, \\[.05cm]
(m)\; f\circ m_A=m_B\circ (f+f), \\[.05cm]
(lm1)\; m_A\circ l^i_{A,A}=\mbox{\bf 1}_A, \\[.05cm]
(lm2)\; m_{A+B}\circ (l^1_{A,B}+l^2_{A,B})=\mbox{\bf 1}_{A+B}.
\end{array}
\]

The isomorphism of ${\cal C}$ and ${\cal C}^{\prime}$ is shown with the
definitions
\[
\begin{array}{ll}
k^1_{A,B}=_{def.} K^1_B\mbox{\bf 1}_A, & l^1_{A,B}=_{def.}L^1_B\mbox{\bf 1}%
_A, \\[.05cm]
k^2_{A,B}=_{def.} K^2_A\mbox{\bf 1}_B, & l^2_{A,B}=_{def.}L^2_A\mbox{\bf 1}%
_B, \\[.05cm]
w_A=_{def.}\langle\mbox{\bf 1}_A,\mbox{\bf 1}_A\rangle, & m_A=_{def.}[%
\mbox{\bf 1}_A,\mbox{\bf 1}_A], \\[.05cm]
f\times g=_{def.}\langle K^1_C f, K^2_A g\rangle, & f+g=_{def.} [L^1_D f,
L^2_B g],
\end{array}
\]
\[
\begin{array}{ll}
K^1_B f=_{def.} f\circ k^1_{A,B}, & L^1_B f=_{def.} l^1_{A,B}\circ f, \\%
[.05cm]
K^2_A f=_{def.} f\circ k^2_{A,B}, & L^2_A f=_{def.} l^2_{A,B}\circ f, \\%
[.05cm]
\langle f, g\rangle =_{def.} (f\times g)\circ w_C, & [f,g]=_{def.}m_C\circ
(f+g).
\end{array}
\]

The free $\times\!,\!+$-{\em categories} ${\cal C}_{\times,+}$ and ${\cal C}%
_{\times,+}^{\prime}$ generated by ${\cal L}$ have as objects formulae of $%
{\cal P}_{\times,+}$. In that case, the terms $k_A$ and $l_A$ are missing,
and the associated equations $(k)$ and $(l)$ are omitted. The remaining
equations are as in ${\cal C}$ and ${\cal C}^{\prime}$. The categories $%
{\cal C}_{\times,+}$ and ${\cal C}_{\times,+}^{\prime}$ are also isomorphic.
We obtain similarly the free cartesian category ${\cal C_{\times,%
\mbox{\scriptsize{\rm I}}}}$ generated by ${\cal L}$, whose objects are from
${\cal P}_{\times, \mbox{\scriptsize{\rm I}}}$, and the isomorphic category $%
{\cal C}_{\times,\mbox{\scriptsize{\rm I}}}^{\prime}$. In that case, we omit
all the terms, operations and equations tied to $+$ and O. For the free
category with binary product ${\cal C}_{\times}$ generated by ${\cal L}$,
and the isomorphic category ${\cal C}_{\times}^{\prime}$, we omit I, $k_A$
and $(k)$ from ${\cal C}_{\times,\mbox{\scriptsize{\rm I}}}$ and ${\cal C}%
_{\times,\mbox{\scriptsize{\rm I}}}^{\prime}$ respectively.

The categories ${\cal C}^+_{\times,\mbox{\scriptsize{\rm I}}}$ and ${\cal C}%
^+_{\times}$ are obtained by extending ${\cal C}_{\times,\mbox{\scriptsize{%
\rm I}}}$ and ${\cal C}_{\times}$ respectively with the connective +, the
operation on terms +, and the equations $(+1)$ and $(+2)$. We obtain ${{\cal %
C}^+_{\times,\mbox{\scriptsize{\rm I}}}}^{\prime}$ and ${{\cal C}^+_{\times}}%
^{\prime}$, which are isomorphic to ${\cal C}^+_{\times,\mbox{\scriptsize{%
\rm I}}}$ and ${\cal C}^+_{\times}$ respectively, by extending in the same
manner ${\cal C}_{\times,\mbox{\scriptsize{\rm I}}}^{\prime}$ and ${\cal C}%
_{\times}^{\prime}$ respectively.

The categories ${\cal C}_{+,\mbox{\scriptsize{\rm O}}}$ and ${\cal C}_{+}$
are isomorphic to ${\cal C}^{op}_{\times,\mbox{\scriptsize{\rm I}}}$ and $%
{\cal C}^{op}_{\times}$ respectively. We obtain the categories ${\cal C}%
^\times_{+,\mbox{\scriptsize{\rm O}}}$ and ${\cal C}^\times_{+}$, which are
isomorphic to ${{\cal C}^+_{\times,\mbox{\scriptsize{\rm I}}}}^{op}$ and ${%
{\cal C}^+_{\times}}^{op}$ respectively, by extending ${\cal C}_{+,%
\mbox{\scriptsize{\rm O}}}$ and ${\cal C}_{+}$ with the connective $\times$,
the operation on terms $\times$, and the equations $(\times 1)$ and $(\times
2)$. We have also isomorphic primed versions of ${\cal C}^\times_{+,%
\mbox{\scriptsize{\rm O}}}$ and ${\cal C}^\times_{+}$.

Up to a certain point in our exposition (noted in Section 4), we shall
distinguish ${\cal C}$ from ${\cal C}^{\prime}$, and analogously for the
other categories derived from ${\cal C}$, which we have now introduced. We
do that until the statements of our auxiliary results are tied to the
nonprimed or to the primed version of the category in question. Once the
necessity for this distinction ceases, we refer to both of these isomorphic
categories by the nonprimed name.

\section{Cut Elimination}

We prove the following theorem for ${\cal C}$.\\[.3cm]
C{\footnotesize UT} E{\footnotesize LIMINATION}.\quad {\em Every term is
equal to a composition-free term.}\\[.2cm]
{P{\footnotesize ROOF}.}\hspace{1em} \nopagebreak
Take a subterm $g\circ f$ of a term such that both $f$ and $g$ are
composition-free. We call such a term a {\em topmost cut}. We show that $%
g\circ f$ is either equal to a composition-free term, or it is equal to a
term all of whose compositions occur in topmost cuts of strictly smaller
length than the length of $g\circ f$. The possibility of eliminating
compositions in topmost cuts, and hence every composition, follows by
induction on the length of topmost cuts.

The cases where $f$ or $g$ is $\mbox{\bf 1}_A$, or $f$ is $l_A$, or $g$ is $%
k_A$, are taken care of by $(cat\: 1)$, $(l)$ and $(k)$. The cases where $f$
is $K^i_A f^{\prime}$ or $g$ is $L^i_A g^{\prime}$ are taken care of by $%
(K1) $ and $(L1)$. And the cases where $f$ is $[f_1, f_2]$ or $g$ is $%
\langle g_1,g_2\rangle$ are taken care of by $(L3)$ and $(K3)$.

The following cases remain. If $f$ is $k_A$, then $g$ is of a form covered
by cases we dealt with above.

If $f$ is $\langle f_1,f_2\rangle$, then $g$ is either of a form covered by
cases above, or $g$ is $K^i_A g^{\prime}$, in which case we apply $(K2)$.

If $f$ is $L^i_A f^{\prime}$, then $g$ is either of a form covered by cases
above, or $g$ is $[g_1,g_2]$, in which case we apply $(L2)$. This covers all
possible cases. \hfill $\Box$ \\[0.2cm]

In this proof we have used all the equations assumed for ${\cal C}$ except $%
(cat\: 2)$, $(K4)$ and $(L4)$.

A portion of this proof suffices to demonstrate Cut Elimination for ${\cal C}%
_{\times,+}$, ${\cal C}_{\times,\mbox{\scriptsize{\rm I}}}$ and ${\cal C}%
_{\times}$. By duality, we also have Cut Elimination for ${\cal C}_{+,%
\mbox{\scriptsize{\rm O}}}$ and ${\cal C}_{+}$. To demonstrate Cut
Elimination for ${\cal C}_{\times,\mbox{\scriptsize{\rm I}}}^+$ we have to
consider the following additional cases.

If $f$ is $k_A$ or $\langle f_1,f_2\rangle$, then $g$ cannot be of the form $%
g_1+g_2$. If $f$ is $f_1+f_2$, and $g$ is not of a form already covered by
cases in the proof above, then $g$ is of the form $g_1+g_2$, in which case
we apply $(+2)$. This covers all possible cases. A portion of this proof
suffices to demonstrate Cut Elimination for ${\cal C}^+_{\times}$. By
duality, we also obtain Cut Elimination for ${\cal C}^\times_{+,%
\mbox{\scriptsize{\rm O}}}$ and ${\cal C}^\times_{+}$.

A composition-free term of ${\cal C}^+_{\times,\mbox{\scriptsize{\rm I}}}$
is reduced to normal form with the following reductions:
\[
\begin{array}{cc}
redexes {\mbox{\hspace{1cm}}} & {\mbox{\hspace{1cm}}} contracta \\[.2cm]
\mbox{\bf 1}_{A\times B} {\mbox{\hspace{1cm}}} & {\mbox{\hspace{1cm}}}
\langle K^1_B \mbox{\bf 1}_A, K^2_A \mbox{\bf 1}_B \rangle \\[.1cm]
\mbox{\bf 1}_{A+B} {\mbox{\hspace{1cm}}} & {\mbox{\hspace{1cm}}} \mbox{\bf 1}%
_A+\mbox{\bf 1}_B \\[.1cm]
\mbox{\bf 1}_{\mbox{\scriptsize{\rm I}}} {\mbox{\hspace{1cm}}} & {%
\mbox{\hspace{1cm}}} k_{\mbox{\scriptsize{\rm I}}} \\[.1cm]
K^i_A \langle f,g\rangle {\mbox{\hspace{1cm}}} & {\mbox{\hspace{1cm}}}
\langle K^i_A f,K^i_A g\rangle \\[.1cm]
K^1_B k_A {\mbox{\hspace{1cm}}} & {\mbox{\hspace{1cm}}} k_{A\times B} \\%
[.1cm]
K^2_A k_B {\mbox{\hspace{1cm}}} & {\mbox{\hspace{1cm}}} k_{A\times B}
\end{array}
\]

These reductions are strongly normalizing. To show that, let $n_1$ be the
number of connectives $\times$, $+$ and I in the indices of identities, and
let $n_2$ be the number of pairs of brackets $\langle,\rangle$ and $k$ terms
within the scope of an operation $K^i_C$ (not necessarily the immediate
scope). Let the degree of a term be $(n_1,n_2)$, and let these degrees be
lexicographically ordered. Then every reduction decreases the degree.

The reduction from $\mbox{\bf 1}_{\mbox{\scriptsize{\rm I}}}$ to $k_{%
\mbox{\scriptsize{\rm I}}}$ and the last two reductions enable us to reduce
every term different from $k_A$ of type $A\rightarrow \mbox{\rm I}$ to $k_A$%
. We disregard these three reductions to reduce to normal form
composition-free terms of ${\cal C}^+_\times$.

All the reductions above are covered by equations of ${\cal C}^+_{\times,%
\mbox{\scriptsize{\rm I}}}$. For the fourth reduction we have the following
derivation in ${\cal C}^+_\times$:
\[
\begin{array}{rcl}
K^i_A \langle f,g\rangle & = & \langle f,g\rangle \circ K^i_A \mbox{\bf 1}%
_B, \quad {\mbox{\rm by $(cat\: 1)$ and $(K1)$}} \\[.1cm]
& = & \langle K^i_A f, K^i_A g\rangle, \quad {%
\mbox{\rm by $(K3)$, $(K1)$
and $(cat\: 1)$.}}
\end{array}
\]

\section{Coherence}

We shall now define a graphical category ${\cal G}$ into which ${\cal C}$
can be mapped. The objects of ${\cal G}$ are finite ordinals. An arrow $%
f:n\rightarrow m$ of ${\cal G}$ will be a binary relation from $n$ to $m$,
i.e. a subset of $n\times m$ with domain $n$ and codomain $m$. The identity $%
\mbox{\bf 1}_n:n\rightarrow n$ of ${\cal G}$ is the identity relation on $n$%
, and composition of arrows is composition of relations.

For an object $A$ of ${\cal C}$, let $|A|$ be the number of occurrences of
propositional letters in $A$. For example, $|(p\times(q+p))+(\mbox{\rm I}%
\times p)|$ is 4.

We now define a functor $G$ from ${\cal C}^{\prime}$ to ${\cal G}$ such that
$G(A)=|A|$. It is clear that $G(A\times B)=G(A+B)=|A|+|B|$. We define $G$ on
arrows inductively:
\[
\begin{array}{lcl}
G(\mbox{\bf 1}_A) & = & \{(x,x):x\in |A|\} = \mbox{\bf 1}_{|A|}, \\[.3cm]
G(k^1_{A,B}) & = & \{(x,x):x\in |A|\}, \\[.05cm]
G(k^2_{A,B}) & = & \{(x+|A|,x):x\in |B|\}, \\[.05cm]
G(w_A) & = & \{(x,x):x\in|A|\} \cup \{(x,x+|A|):x\in|A|\}, \\[.05cm]
G(k_A) & = & \emptyset,
\end{array}
\]
\[
\begin{array}{lcl}
G(l^1_{A,B}) & = & \{(x,x):x\in |A|\}, \\[.05cm]
G(l^2_{A,B}) & = & \{(x,x+|A|):x\in |B|\}, \\[.05cm]
G(m_A) & = & \{(x,x) :x\in|A|\} \cup \{(x+|A|,x):x\in |A|\}, \\[.05cm]
G(l_A) & = & \emptyset, \\[.3cm]
G(g\circ f) & = & G(g) \circ G(f),
\end{array}
\]
and for $f:A\rightarrow B$ and $g:C\rightarrow D$,
\[
G(f\times g)=G(f+g)=G(f) \cup \{(x+|A|,y+|B|):(x,y)\in G(g)\}.
\]

Though $G(\mbox{\bf 1}_A)$, $G(k^1_{A,B})$ and $G(l^1_{A,B})$ are the same
as sets of ordered pairs, in general they have different domains and
codomains, the first being a subset of $|A|\times |A|$, the second a subset
of $(|A|+|B|)\times |A|$, and the third a subset of $|A|\times(|A|+|B|)$. We
have an analogous situation in some other cases.

It is easy to draw $G(f)$ diagrammatically. For example, for $%
G(m_{p+q}\circ(l^1_{p,q}+l^2_{p,q}))$ we have

\begin{center}
\begin{picture}(380,80)
\put(50,10){\makebox(0,0){$p$}}
\put(60,10){\makebox(0,0){$+$}}
\put(70,10){\makebox(0,0){$q$}}

\put(10,40){\makebox(0,0){$($}}
\put(20,40){\makebox(0,0){$p$}}
\put(30,40){\makebox(0,0){$+$}}
\put(40,40){\makebox(0,0){$q$}}
\put(50,40){\makebox(0,0){$)$}}
\put(60,40){\makebox(0,0){$+$}}
\put(70,40){\makebox(0,0){$($}}
\put(80,40){\makebox(0,0){$p$}}
\put(90,40){\makebox(0,0){$+$}}
\put(100,40){\makebox(0,0){$q$}}
\put(110,40){\makebox(0,0){$)$}}

\put(50,70){\makebox(0,0){$p$}}
\put(60,70){\makebox(0,0){$+$}}
\put(70,70){\makebox(0,0){$q$}}

\put(50,65){\line(-3,-2){30}}
\put(70,65){\line(3,-2){30}}

\put(20,35){\line(3,-2){28}}
\put(80,35){\line(-3,-2){28}}
\put(40,35){\line(3,-2){28}}
\put(100,35){\line(-3,-2){28}}

\put(150,55){\makebox(0,0){$l^1_{p,q}+l^2_{p,q}$}}
\put(150,25){\makebox(0,0){$m_{p+q}$}}
\put(220,40){\makebox(0,0){which is equal to}}

\put(280,20){\makebox(0,0){$p$}}
\put(290,20){\makebox(0,0){$+$}}
\put(300,20){\makebox(0,0){$q$}}

\put(280,60){\makebox(0,0){$p$}}
\put(290,60){\makebox(0,0){$+$}}
\put(300,60){\makebox(0,0){$q$}}

\put(280,55){\line(0,-1){30}}
\put(300,55){\line(0,-1){30}}

\put(325,40){\makebox(0,0){${\bf 1}_{p+q}$}}
\end{picture}
\end{center}

It is also easy to check that $G$ is a functor from ${\cal C}^{\prime}$ to $%
{\cal G}$. We show by induction on the length of derivation that if $f=g$ in
${\cal C}^{\prime}$, then $G(f)=G(g)$ in ${\cal G}$. (Of course, $G$
preserves identities and composition.) Since, ${\cal C}^{\prime}$ and ${\cal %
C}$ are isomorphic we also have a functor from ${\cal C}$ to ${\cal G}$.

For the bicartesian structure of ${\cal G}$ we have that the operations $%
\times$ and $+$ on objects are both addition of ordinals, the operations $%
\times$ and $+$ on arrows coincide and are defined by the clauses for $%
G(f\times g)$ and $G(f+g)$, and the terminal and the initial object also
coincide: they are both the ordinal zero. The category ${\cal G}$ has zero
arrows, namely, the arrows
\[
\begin{array}{ccccccc}
n & \stackrel{\emptyset}{\rightarrow} & \emptyset & \stackrel{\emptyset}{%
\rightarrow} & \emptyset & \stackrel{\emptyset}{\rightarrow} & m \\
\parallel &  & \parallel &  & \parallel &  & \parallel \\
G(A) &  & G(\mbox{\rm I}) &  & G(\mbox{\rm O}) &  & G(B)
\end{array}
\]
which composed with any other arrow give another zero arrow. It is easy to
see that the bicartesian category ${\cal G}$ is a linear category in the
sense of \cite{law97} (see p. 279). The functor $G$ from ${\cal C}$ to $%
{\cal G}$ is not just a functor, but a bicartesian functor; namely, a
functor that preserves the bicartesian structure of ${\cal C}$.

We also have functors defined analogously to $G$, which we call $G$ too,
from ${\cal C}_{\times,+}$, ${\cal C}^+_{\times,\mbox{\scriptsize{\rm I}}}$
and ${\cal C}^+_{\times}$ to ${\cal G}$. These functors, which are defined
officially for the primed versions of these categories, are obtained from
the definition of $G$ above by just rejecting clauses that are no longer
applicable. For these last three functors we shall show that they are
faithful. By duality, we also have faithful functors $G$ from ${\cal C}%
^\times_{+, \mbox{\scriptsize{\rm O}}}$ and ${\cal C}^\times_{+}$ to ${\cal G%
}$.

That an analogously defined functor $G$ exists from ${\cal C}_{\times, %
\mbox{\scriptsize{\rm I}}}$ to ${\cal G}$, and is faithful, has been
announced in \cite{kel72} (p. 129) and proved in \cite{min80} (Theorem 2.2),
\cite{tro96} (Theorem 8.2.3, p. 207), \cite{pet98} and \cite{dos00}. The
functor $G$ from ${\cal C}_{\times,\mbox{\scriptsize{\rm I}}}$ maps ${\cal C}%
_{\times,\mbox{\scriptsize{\rm I}}}$ into the subcategory of ${\cal G}$
whose arrows are relations converse to functions; in other words, $G$ maps $%
{\cal C}_{+,\mbox{\scriptsize{\rm O}}}$ into the subcategory of ${\cal G}$
whose arrows are functions.

It is clear that the functor $G$ from ${\cal C}$ to ${\cal G}$ is not full,
since there are no arrows in ${\cal C}$ from I to O. This functor is also
not faithful. The counterexamples that show that are
\[
\begin{array}{l}
G(k^1_{\mbox{\scriptsize{\rm O}},\mbox{\scriptsize{\rm O}}})=G(k^2_{%
\mbox{\scriptsize{\rm O}},\mbox{\scriptsize{\rm O}}})=\emptyset, \\[.05cm]
G(l^1_{\mbox{\scriptsize{\rm I}},\mbox{\scriptsize{\rm I}}})=G(l^2_{%
\mbox{\scriptsize{\rm I}},\mbox{\scriptsize{\rm I}}})=\emptyset,
\end{array}
\]
whereas $k^1_{\mbox{\scriptsize{\rm O}},\mbox{\scriptsize{\rm O}}}=k^2_{%
\mbox{\scriptsize{\rm O}},\mbox{\scriptsize{\rm O}}}$ and $l^1_{%
\mbox{\scriptsize{\rm I}},\mbox{\scriptsize{\rm I}}}=l^2_{%
\mbox{\scriptsize{\rm I}},\mbox{\scriptsize{\rm I}}}$ don't hold in ${\cal C}
$. That these equations don't hold in ${\cal C}$ is demonstrated by the
bicartesian category {\bf Set} of sets with functions, with cartesian
product $\times$, disjoint union $+$, singleton I and empty set O. In {\bf %
Set} we don't have $l^1_{\mbox{\scriptsize{\rm I}},\mbox{\scriptsize{\rm I}}%
}=l^2_{\mbox{\scriptsize{\rm I}},\mbox{\scriptsize{\rm I}}}$, and in the
bicartesian category ${\makebox{\bf Set}}^{op}$ we don't have $k^1_{%
\mbox{\scriptsize{\rm O}},\mbox{\scriptsize{\rm O}}}=k^2_{%
\mbox{\scriptsize{\rm O}},\mbox{\scriptsize{\rm O}}}$. The same
counterexamples show that the functors $G$ from ${\cal C}_{\times,+,%
\mbox{\scriptsize{\rm I}}}$ to ${\cal G}$, and from ${\cal C}_{\times,+,%
\mbox{\scriptsize{\rm O}}}$ to ${\cal G}$, are not faithful.

To prove the faithfulness of $G$ from ${\cal C}^+_{\times,%
\mbox{\scriptsize{\rm I}}}$ we need the following lemma.\\[.3cm]
{L{\footnotesize EMMA} 4.1.}\hspace{1em} {\em If} $f,g:A\rightarrow B$ {\em %
are composition-free terms of} ${\cal C}^+_{\times,\mbox{\scriptsize{\rm I}}%
} $ {\em in normal form and} $G(f)=G(g)$ {\em in} ${\cal G}$, {\em then} $f$
{\em and} $g$ {\em are the same term.} \\[.2cm]
{P{\footnotesize ROOF}.}\hspace{1em} \nopagebreak We proceed by induction on
the length of $f$. If $f$ is $\mbox{\bf 1}_p$, then $g$ must be $\mbox{\bf 1}%
_p$, and if $f$ is $k_A$, then $g$ must be $k_A$. If $f$ is $%
K^{i_1}_{A_1}\ldots K^{i_n}_{A_n}f^{\prime}$ for $n\geq 1$, and $f^{\prime}$
is either $\mbox{\bf 1}_p$ or $f_1^{\prime}+f_2^{\prime}$, then with the
help of $G(f)=G(g)$ we conclude that $g$ too must be of the form $%
K^{i_1}_{A_1}\ldots K^{i_n}_{A_n}g^{\prime}$ for $g^{\prime}$ either $%
\mbox{\bf 1}_p$ or $g_1^{\prime}+g_2^{\prime}$. In the latter case, we apply
the induction hypothesis. We apply the induction hypothesis also when $f$ is
$\langle f_1,f_2\rangle$ or $f_1+f_2$. \hfill $\Box$ \\[0.2cm]

As a corollary we obtain the following proposition.\\[.3cm]
U{\footnotesize NIQUENESS OF} C{\footnotesize OMPOSITION}-F{\footnotesize REE%
} N{\footnotesize ORMAL} F{\footnotesize ORM FOR} ${\cal C}^+_{\times,%
\mbox{\scriptsize{\rm I}}}$.\quad {\em If} $f=g$ {\em in} ${\cal C}%
^+_{\times,\mbox{\scriptsize{\rm I}}}$ {\em for} $f$ {\em and} $g$ {\em in
composition-free normal form, then} $f$ {\em and} $g$ {\em are the same term}%
.\\[.2cm]
{P{\footnotesize ROOF}.}\hspace{1em} \nopagebreak From $f=g$ it follows that
$f$ and $g$ are of the same type and that $G(f)=G(g)$ in ${\cal G}$. Then we
apply Lemma 4.1. \hfill $\Box$ \\[0.2cm]

Note that we have established this uniqueness without appealing to the
Church-Rosser property for our reductions. Another corollary of Lemma 4.1 is
that $(cat\: 2)$ can be derived from the remaining equations of ${\cal C}%
^+_{\times,\mbox{\scriptsize{\rm I}}}$. To derive $h\circ(g\circ f)=(h\circ
g)\circ f$, we reduce both sides to cut-free normal form, which is done
without using $(cat\: 2)$. An analogous lemma, which implies Uniqueness of
Composition-Free Normal Form and derivability of $(cat\: 2)$, can also be
established for ${\cal C}^+_{\times}$, ${\cal C}_{\times,\mbox{\scriptsize{%
\rm I}}}$ and ${\cal C}_{\times}$.

We can now establish the following coherence proposition.\\[.3cm]
F{\footnotesize AITHFULNESS OF} $G$ {\footnotesize FROM} ${\cal C}^+_{\times,%
\mbox{\scriptsize{\rm I}}}$. \quad {\em If} $f,g:A\rightarrow B$ {\em are
terms of} ${\cal C}^+_{\times,\mbox{\scriptsize{\rm I}}}$ {\em and} $%
G(f)=G(g)$ {\em in} ${\cal G}$, {\em then} $f=g$ {\em in} ${\cal C}%
^+_{\times,\mbox{\scriptsize{\rm I}}}$.\\[.2cm]
{P{\footnotesize ROOF}.}\hspace{1em} \nopagebreak Suppose $f,g:A\rightarrow
B $ are terms of ${\cal C}^+_{\times,\mbox{\scriptsize{\rm I}}}$, and $%
f^{\prime}$ and $g^{\prime}$ are the composition-free normal forms of $f$
and $g$ respectively. Then from $G(f)=G(g)$, $G(f)=G(f^{\prime})$ and $%
G(g)=G(g^{\prime})$ we obtain $G(f^{\prime})=G(g^{\prime})$, and therefore,
by Lemma 4.1, it follows that $f^{\prime}$ and $g^{\prime}$ are the same
term. Hence $f=g$ in ${\cal C}^+_{\times,\mbox{\scriptsize{\rm I}}}$. \hfill
$\Box$ \\[0.2cm]

A portion of this proof suffices to demonstrate that the functor $G$ from $%
{\cal C}^+_{\times}$ to ${\cal G}$ is also faithful, and we can also
demonstrate in the same manner the faithfulness of the functor $G$ from $%
{\cal C}_{\times,\mbox{\scriptsize{\rm I}}}$ to ${\cal G}$, or from ${\cal C}%
_{\times}$ to ${\cal G}$, but this is already known, as we noted above. By
duality, we also have the faithfulness of $G$ from ${\cal C}^\times_{+,%
\mbox{\scriptsize{\rm O}}}$, ${\cal C}^\times_{+}$, ${\cal C}_{+,%
\mbox{\scriptsize{\rm O}}}$ and ${\cal C}_{\mbox{\scriptsize{\rm O}}}$. It
remains to demonstrate that the functor $G$ from ${\cal C}_{\times,+}$ to $%
{\cal G}$ is also faithful.

For a term of ${\cal C}^{\prime}$ of the form $f_n\circ \ldots\circ f_1$,
for some $n\geq 1$, where $f_i$ is composition-free we shall say that it is
{\em factorized}. By using $(\times 2)$, $(+2)$ and $(cat\: 1)$ it is easy
to show that every term of ${\cal C}^{\prime}$ is equal to a factorized term
of ${\cal C}^{\prime}$. A subterm $f_i$ in a factorized term $f_n\circ
\ldots\circ f_1$ is called a {\em factor}.

A term of ${\cal C}^{\prime}$ where all the atomic terms are identities will
be called a {\em complex identity}. According to $(\times 1)$, $(+1)$ and $%
(cat\: 1)$, every complex identity is equal to an identity. A factor which
is a complex identity will be called an {\em identity factor}. It is clear
that if $n>1$, we can omit in a factorized term every identity factor, and
obtain a factorized term equal to the original one.

A term of ${\cal C}^{\prime}_{\times,+}$ is called a {\em K-term} iff it is
a term of ${{\cal C}^+_{\times}}^{\prime}$ and it is not a complex identity.
A term of ${\cal C}^{\prime}_{\times,+}$ is called an {\em L-term} iff it is
a term of ${{\cal C}_+^{\times}}^{\prime}$ and it is not a complex identity.
Remember that the terms of ${{\cal C}^+_{\times}}^{\prime}$ have the atomic
terms $\mbox{\bf 1}_A,k^i_{A,B}$ and $w_A$, and the operations on terms $%
\circ, \times$ and $+$; the terms of ${{\cal C}_+^{\times}}^{\prime}$ have
the atomic terms $\mbox{\bf 1}_A, l^i_{A,B}$ and $m_A$, and the same
operations on terms.

A term of ${\cal C}^{\prime}_{\times,+}$ is said to be in {\em K-L normal
form} iff it is of the form $g\circ f:A\rightarrow B$ for $f$ a {\em K}-term
or $\mbox{\bf 1}_A$ and $g$ an {\em L}-term or $\mbox{\bf 1}_B$. Note that
{\em K-L} normal forms are not unique, since $(m_A\times m_A)\circ w_{A+A}$
and $m_{A\times A}\circ(w_A+w_A)$, which are both equal to $w_A\circ m_A$,
are both in {\em K-L} normal form.

We can prove the following proposition.\\[.3cm]
{\em K-L} N{\footnotesize ORMALIZATION}. \quad {\em Every term of} ${\cal C}%
^{\prime}_{\times,+}$ {\em is equal in} ${\cal C}^{\prime}_{\times,+}$ {\em %
to a term of} ${\cal C}^{\prime}_{\times,+}$ {\em in} {\em K-L} {\em normal
form}.\\[.2cm]
{P{\footnotesize ROOF}.}\hspace{1em} \nopagebreak Suppose $f:B\rightarrow C$
is a composition-free {\em K}-term and $g:A\rightarrow B$ is a
composition-free {\em L}-term. We show by induction on the length of $f\circ
g$ that
\[
(*)\quad f\circ g=g^{\prime}\circ f^{\prime}\; {\mbox{\rm or }} f\circ
g=f^{\prime}\; {\mbox{\rm or }} f\circ g=g^{\prime}
\]
for $f^{\prime}$ a composition-free {\em K}-term and $g^{\prime}$ a
composition-free {\em L}-term.

We shall not consider below cases where $g$ is $m_B$, which are easily taken
care of by $(m)$. The following cases remain.

If $f$ is $k^i_{C,E}$ and $g$ is $g_1\times g_2$, then we use $(k^i)$. If $f$
is $w_B$, then we use $(w)$. If $f$ is $f_1\times f_2$ and $g$ is $g_1\times
g_2$, then we use $(\times 2)$, the induction hypothesis, and perhaps $%
(cat\: 1)$.

Finally, if $f$ is $f_1+f_2$, then we have the following cases. If $g$ is $%
l^i_{B_1,B_2}$, then we use $(l^i)$. If $g$ is $g_1+g_2$, then we use $(+2)$%
, the induction hypothesis, and perhaps $(cat\: 1)$. This proves $(*)$.

Every term of ${\cal C}_{\times,+}^{\prime}$ is equal to an identity or to a
factorized term $f_n\circ\ldots\circ f_1$ without identity factors. Every
factor $f_i$ of $f_n\circ\ldots\circ f_1$ is either a {\em K}-term or an
{\em L}-term or, by $(cat\: 1)$, $(\times2)$ and $(+2)$, it is equal to $%
f_i^{\prime\prime}\circ f_i^{\prime}$ where $f_i^{\prime}$ is a {\em K}-term
and $f_i^{\prime\prime}$ is an {\em L}-term. For example, $(k^1_{A,B}\times
l^1_{C,D})+w_E$ is equal to
\[
((\mbox{\bf 1}_A\times l^1_{C,D})+\mbox{\bf 1}_{E\times E})\circ
((k^1_{A,B}\times\mbox{\bf 1}_C)+w_E).
\]
Then it is clear that by applying $(*)$ repeatedly, and by applying perhaps $%
(cat\: 1)$ at the end, we obtain a term in {\em K-L} normal form. \hfill $%
\Box$ \\[0.2cm]

Note that to reduce a term of ${\cal C}_{\times,+}^{\prime}$ to {\em K-L}
normal form we have used in this proof all the equations of ${\cal C}%
_{\times,+}^{\prime}$ except $(kw1)$, $(kw2)$, $(lm1)$ and $(lm2)$.

A term of ${\cal C}^{\prime}$ is called a {\em K-term} iff $l_A$, $l^i_{A,B}$
and $m_A$ don't occur in it and it is not a complex identity. A term of $%
{\cal C}^{\prime}$ is called an {\em L-term} iff $k_A$, $k^i_{A,B}$ and $w_A$
don't occur in it and it is not a complex identity. The definition of {\em %
K-L} normal form is as above. Then we can prove {\em K-L} Normalization for $%
{\cal C}^{\prime}$ too. It is enough to consider in the induction that
establishes $(*)$ in the proof above the additional cases where $f$ is $k_B$
or $g$ is $l_B$, which are easily taken care of by $(k)$ and $(l)$.

From now on we shall make no distinction any more between the categories $%
{\cal C}$ and ${\cal C}^{\prime}$. These categories are isomorphic, and both
will be called ${\cal C}$. When we refer, for example, to the term $%
k^1_{A,B} $ of ${\cal C}$, we refer to the term defined as in Section 2, and
analogously in other cases. We proceed in the same way in making no
distinction between other categories derived from ${\cal C}$ and their
isomorphic primed versions. We can then prove the following lemma.\\[.3cm]
{L{\footnotesize EMMA} 4.2.}\hspace{1em} {\em Let} $f:A_1\times
A_2\rightarrow B$ {\em be a term of} ${\cal C}^+_{\times}$. {\em If for every%
} $(x,y)\in G(f)$ {\em we have} $x\in|A_1|$, {\em then} $f$ {\em is equal in}
${\cal C}^+_{\times}$ {\em to a term of the form} $f^{\prime}\circ
k^1_{A_1,A_2}$, {\em and if for every} $(x,y)\in G(f)$ {\em we have} $%
x-|A_1|\in|A_2|$, {\em then} $f$ {\em is equal in} ${\cal C}^+_{\times}$
{\em to a term of the form} $f^{\prime}\circ k^2_{A_1,A_2}$.\\[.2cm]
{P{\footnotesize ROOF}.}\hspace{1em} \nopagebreak We proceed by induction on
the length of $B$. If $B$ is a propositional letter or $B_1+B_2$, then by
Cut Elimination $f$ must be equal to a term of the form $f^{\prime}\circ
k^i_{A_1,A_2}$. The condition on $G(f)$ dictates whether $i$ here is 1 or 2.

If $B$ is $B_1\times B_2$, and for every $(x,y)\in G(f)$ we have $x\in |A_1|$%
, then for $k^i_{B_1,B_2}\circ f:A_1\times A_2\rightarrow B_i$, for every $%
(x,z)\in G(k^i_{B_1,B_2}\circ f)$ we have $x\in |A_1|$. So, by the induction
hypothesis,
\[
k^i_{B_1,B_2}\circ f=f_i\circ k^1_{A_1,A_2}.
\]
Hence
\[
\begin{array}{rcl}
f & = & \langle k^1_{B_1,B_2}\circ f, k^2_{B_1,B_2}\circ f\rangle \\
& = & \langle f_1,f_2\rangle\circ k^1_{A_1,A_2}.
\end{array}
\]

We reason analogously if for every $(x,y)\in G(f)$ we have $x-|A_1|\in |A_2|$%
. \hfill $\Box$ \\[0.2cm]

We can prove analogously the following dual lemma.\\[.3cm]
{L{\footnotesize EMMA} 4.3.}\hspace{1em} {\em Let} $f:A\rightarrow B_1+B_2$
{\em be a term of} ${\cal C}^{\times}_{+}$. {\em If for every} $(x,y)\in
G(f) $ {\em we have} $y\in |B_1|$, {\em then} $g$ {\em is equal in} ${\cal C}%
^{\times}_{+}$ {\em to a term of the form} $l^1_{B_1,B_2}\circ g^{\prime}$,
{\em and if for every} $(x,y)\in G(f)$ {\em we have} $y-|B_1|\in |B_2|$,
{\em then} $g$ {\em is equal in} ${\cal C}^{\times}_{+}$ {\em to a term of
the form} $l^2_{B_1,B_2}\circ g^{\prime}$.\\[.2cm]

We shall next prove the following coherence proposition.\\[.3cm]
F{\footnotesize AITHFULNESS OF} $G$ {\footnotesize FROM} ${\cal C}%
_{\times,+} $. \quad {\em If} $f,g:A\rightarrow B$ {\em are terms of} ${\cal %
C}_{\times,+}$ {\em and} $G(f)=G(g)$ {\em in} ${\cal G}$, {\em then} $f=g$
{\em in} ${\cal C}_{\times,+}$. \\[.2cm]
{P{\footnotesize ROOF}.}\hspace{1em} \nopagebreak We proceed by induction on
the length of $A$ and $B$. In the basis of this induction, when both $A$ and
$B$ are propositional letters, we conclude by Cut Elimination that $f$ and $%
g $ exist iff $A$ and $B$ are the same propositional letter $p$, and $f=g=%
\mbox{\bf 1}_p$ in ${\cal C}_{\times,+}$. (We could conclude the same thing
by interpreting ${\cal C}_{\times,+}$ in conjunctive-disjunctive logic.)
Note that we didn't need here the assumption $G(f)=G(g)$.

If $A$ is $A_1+A_2$, then $f\circ l^1_{A_1,A_2}$ and $g\circ l^1_{A_1,A_2}$
are of type $A_1\rightarrow B$, while $f\circ l^2_{A_1,A_2}$ and $g\circ
l^2_{A_1,A_2}$ are of type $A_2\rightarrow B$. We also have
\[
\begin{array}{rcl}
G(f\circ l^i_{A_1,A_2}) & = & G(f)\circ G(l^i_{A_1,A_2}) \\
& = & G(g)\circ G(l^i_{A_1,A_2}) \\
& = & G(g\circ l^i_{A_1,A_2}),
\end{array}
\]
whence, by the induction hypothesis,
\[
f\circ l^i_{A_1,A_2}=g\circ l^i_{A_1,A_2}
\]
in ${\cal C}_{\times,+}$. Then we infer that
\[
[f\circ l^1_{A_1,A_2},f\circ l^2_{A_1,A_2}]=[g\circ l^1_{A_1,A_2},g\circ
l^2_{A_1,A_2}],
\]
from which it follows that $f=g$ in ${\cal C}_{\times,+}$. We proceed
analogously if $B$ is $B_1\times B_2$.

Suppose now $A$ is $A_1\times A_2$ or a propositional letter, and $B$ is $%
B_1+B_2$ or a propositional letter, but $A$ and $B$ are not both
propositional letters. Then, by Cut Elimination, $f$ is equal either in $%
{\cal C}_{\times,+}$ to a term of the form $f^{\prime}\circ k^i_{A_1,A_2}$,
or to a term of the form $l^i_{B_1,B_2}\circ f^{\prime}$. Suppose $%
f=f^{\prime}\circ k^1_{A_1,A_2}$. Then for every $(x,y)\in G(f)$ we have $%
x\in |A_1|$. (We reason analogously when $f=f^{\prime}\circ k^2_{A_1,A_2}$.)

By {\em K-L} normalization, $g=g_2\circ g_1$ in ${\cal C}_{\times,+}$ for $%
g_1:A_1\times A_2\rightarrow C$ a term of ${\cal C}^+_{\times}$ and $g_2$ a
term of ${\cal C}_+^{\times}$. Since $g_2:C\rightarrow B$ is a term of $%
{\cal C}_+^{\times}$, and hence $K^i$ (that is, $k^i$) does not occur in it,
for every $z\in |C|$ we have a $y\in|B|$ such that $(z,y)\in G(g_2)$. If for
some $(x,z)\in G(g_1)$ we had $x\not\in |A_1|$, then for some $(x,y)\in
G(g_2\circ g_1)$ we would have $x\not\in |A_1|$, but this is impossible
since $G(g_2\circ g_1)=G(g)=G(f)$. So for every $(x,z)\in G(g_1)$ we have $%
x\in|A_1|$. Then, by Lemma 4.2, $g_1=g_1^{\prime}\circ k^1_{A_1,A_2}$ in $%
{\cal C}^+_{\times}$, and hence in ${\cal C}_{\times,+}$ too. Therefore, $%
g=g_2\circ g_1^{\prime}\circ k^1_{A_1,A_2}$ in ${\cal C}_{\times,+}$.

Because of the particular form of $G(k^1_{A_1,A_2})$, we can infer from $%
G(f)=G(g)$ that $G(f^{\prime})=G(g_2\circ g_1^{\prime})$, but since $%
f^{\prime}$ and $g_2\circ g_1^{\prime}$ are of type $A_1\rightarrow B$, by
the induction hypothesis we have $f^{\prime}=g_2\circ g_1^{\prime}$ in $%
{\cal C}_{\times,+}$, and hence $f=g$. When $f=l^i_{B_1,B_2}\circ f^{\prime}$%
, we reason analogously and apply Lemma 4.3. \hfill $\Box$ \\[0.2cm]

With the help of the Faithfulness of $G$ from ${\cal C}_{\times,+}$ it is
easy to establish, for example, that in ${\cal C}_{\times,+}$
\[
\langle [f_1,f_2],[g_1,g_2]\rangle = [\langle f_1,g_1\rangle ,\langle
f_2,g_2\rangle ],
\]
for which we have the diagrams

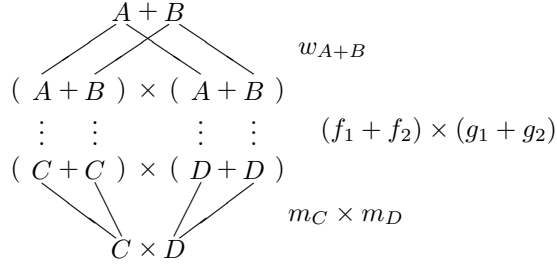
\begin{center}
\begin{picture}(200,110)
\put(50,10){\makebox(0,0){$C$}}
\put(60,10){\makebox(0,0){$\times$}}
\put(70,10){\makebox(0,0){$D$}}

\put(10,40){\makebox(0,0){$($}}
\put(20,40){\makebox(0,0){$C$}}
\put(30,40){\makebox(0,0){$+$}}
\put(40,40){\makebox(0,0){$C$}}
\put(50,40){\makebox(0,0){$)$}}
\put(60,40){\makebox(0,0){$\times$}}
\put(70,40){\makebox(0,0){$($}}
\put(80,40){\makebox(0,0){$D$}}
\put(90,40){\makebox(0,0){$+$}}
\put(100,40){\makebox(0,0){$D$}}
\put(110,40){\makebox(0,0){$)$}}

\put(10,70){\makebox(0,0){$($}}
\put(20,70){\makebox(0,0){$A$}}
\put(30,70){\makebox(0,0){$+$}}
\put(40,70){\makebox(0,0){$B$}}
\put(50,70){\makebox(0,0){$)$}}
\put(60,70){\makebox(0,0){$\times$}}
\put(70,70){\makebox(0,0){$($}}
\put(80,70){\makebox(0,0){$A$}}
\put(90,70){\makebox(0,0){$+$}}
\put(100,70){\makebox(0,0){$B$}}
\put(110,70){\makebox(0,0){$)$}}

\put(50,100){\makebox(0,0){$A$}}
\put(60,100){\makebox(0,0){$+$}}
\put(70,100){\makebox(0,0){$B$}}

\put(20,35){\line(4,-3){28}}
\put(80,35){\line(-1,-2){10}}
\put(40,35){\line(1,-2){10}}
\put(100,35){\line(-4,-3){28}}

\put(20,57){\makebox(0,0){$\vdots$}}
\put(40,57){\makebox(0,0){$\vdots$}}
\put(80,57){\makebox(0,0){$\vdots$}}
\put(100,57){\makebox(0,0){$\vdots$}}

\put(48,95){\line(-3,-2){28}}
\put(52,95){\line(3,-2){28}}
\put(68,95){\line(-3,-2){28}}
\put(72,95){\line(3,-2){28}}

\put(135,23){\makebox(0,0){$m_C\times m_D$}}
\put(170,55){\makebox(0,0){$(f_1+f_2)\times (g_1+g_2)$}}
\put(130,85){\makebox(0,0){$w_{A+B}$}}
\end{picture}

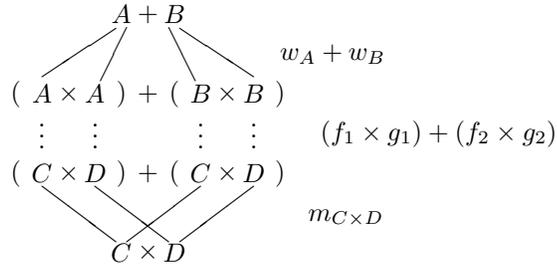
\begin{picture}(200,110)
\put(50,10){\makebox(0,0){$C$}}
\put(60,10){\makebox(0,0){$\times$}}
\put(70,10){\makebox(0,0){$D$}}

\put(10,40){\makebox(0,0){$($}}
\put(20,40){\makebox(0,0){$C$}}
\put(30,40){\makebox(0,0){$\times$}}
\put(40,40){\makebox(0,0){$D$}}
\put(50,40){\makebox(0,0){$)$}}
\put(60,40){\makebox(0,0){$+$}}
\put(70,40){\makebox(0,0){$($}}
\put(80,40){\makebox(0,0){$C$}}
\put(90,40){\makebox(0,0){$\times$}}
\put(100,40){\makebox(0,0){$D$}}
\put(110,40){\makebox(0,0){$)$}}

\put(10,70){\makebox(0,0){$($}}
\put(20,70){\makebox(0,0){$A$}}
\put(30,70){\makebox(0,0){$\times$}}
\put(40,70){\makebox(0,0){$A$}}
\put(50,70){\makebox(0,0){$)$}}
\put(60,70){\makebox(0,0){$+$}}
\put(70,70){\makebox(0,0){$($}}
\put(80,70){\makebox(0,0){$B$}}
\put(90,70){\makebox(0,0){$\times$}}
\put(100,70){\makebox(0,0){$B$}}
\put(110,70){\makebox(0,0){$)$}}

\put(50,100){\makebox(0,0){$A$}}
\put(60,100){\makebox(0,0){$+$}}
\put(70,100){\makebox(0,0){$B$}}

\put(20,35){\line(4,-3){28}}
\put(80,35){\line(-4,-3){28}}
\put(40,35){\line(4,-3){28}}
\put(100,35){\line(-4,-3){28}}

\put(20,57){\makebox(0,0){$\vdots$}}
\put(40,57){\makebox(0,0){$\vdots$}}
\put(80,57){\makebox(0,0){$\vdots$}}
\put(100,57){\makebox(0,0){$\vdots$}}

\put(48,95){\line(-3,-2){28}}
\put(52,95){\line(-1,-2){10}}
\put(68,95){\line(1,-2){10}}
\put(72,95){\line(3,-2){28}}

\put(135,23){\makebox(0,0){$m_{C\times D}$}}
\put(170,55){\makebox(0,0){$(f_1\times g_1)+ (f_2\times g_2)$}}
\put(130,85){\makebox(0,0){$w_A+w_B$}}
\end{picture}
\end{center}

\noindent or that in ${\cal C}_{\times,+}$
\[
\begin{array}{c}
((k^1_{A,B}+k^1_{C,D})\times(k^2_{A,B}+k^2_{C,D}))\circ w_{(A\times
B)+(C\times D)}= \\
{\mbox{\hspace{3em}}}= m_{(A+C)\times (B+D)}\circ ((l^1_{A,C}\times
l^1_{B,D})+(l^2_{A,C}\times l^2_{B,D})),
\end{array}
\]
for which we have the diagrams

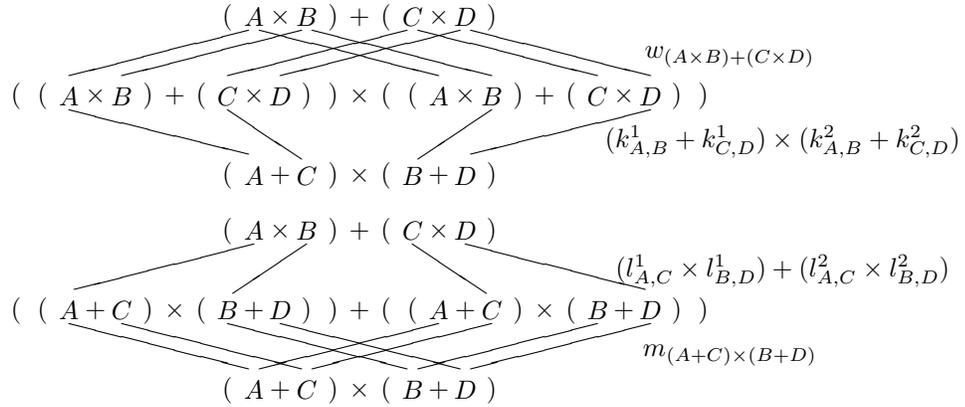
\begin{center}
\begin{picture}(350,80)
\put(90,10){\makebox(0,0){$($}}
\put(100,10){\makebox(0,0){$A$}}
\put(110,10){\makebox(0,0){$+$}}
\put(120,10){\makebox(0,0){$C$}}
\put(130,10){\makebox(0,0){$)$}}
\put(140,10){\makebox(0,0){$\times$}}
\put(150,10){\makebox(0,0){$($}}
\put(160,10){\makebox(0,0){$B$}}
\put(170,10){\makebox(0,0){$+$}}
\put(180,10){\makebox(0,0){$D$}}
\put(190,10){\makebox(0,0){$)$}}

\put(10,40){\makebox(0,0){$($}}
\put(20,40){\makebox(0,0){$($}}
\put(30,40){\makebox(0,0){$A$}}
\put(40,40){\makebox(0,0){$\times$}}
\put(50,40){\makebox(0,0){$B$}}
\put(60,40){\makebox(0,0){$)$}}
\put(70,40){\makebox(0,0){$+$}}
\put(80,40){\makebox(0,0){$($}}
\put(90,40){\makebox(0,0){$C$}}
\put(100,40){\makebox(0,0){$\times$}}
\put(110,40){\makebox(0,0){$D$}}
\put(120,40){\makebox(0,0){$)$}}
\put(130,40){\makebox(0,0){$)$}}
\put(140,40){\makebox(0,0){$\times$}}
\put(150,40){\makebox(0,0){$($}}
\put(160,40){\makebox(0,0){$($}}
\put(170,40){\makebox(0,0){$A$}}
\put(180,40){\makebox(0,0){$\times$}}
\put(190,40){\makebox(0,0){$B$}}
\put(200,40){\makebox(0,0){$)$}}
\put(210,40){\makebox(0,0){$+$}}
\put(220,40){\makebox(0,0){$($}}
\put(230,40){\makebox(0,0){$C$}}
\put(240,40){\makebox(0,0){$\times$}}
\put(250,40){\makebox(0,0){$D$}}
\put(260,40){\makebox(0,0){$)$}}
\put(270,40){\makebox(0,0){$)$}}

\put(90,70){\makebox(0,0){$($}}
\put(100,70){\makebox(0,0){$A$}}
\put(110,70){\makebox(0,0){$\times$}}
\put(120,70){\makebox(0,0){$B$}}
\put(130,70){\makebox(0,0){$)$}}
\put(140,70){\makebox(0,0){$+$}}
\put(150,70){\makebox(0,0){$($}}
\put(160,70){\makebox(0,0){$C$}}
\put(170,70){\makebox(0,0){$\times$}}
\put(180,70){\makebox(0,0){$D$}}
\put(190,70){\makebox(0,0){$)$}}

\put(30,35){\line(4,-1){68}}
\put(90,35){\line(3,-2){28}}
\put(190,35){\line(-3,-2){28}}
\put(250,35){\line(-4,-1){68}}

\put(98,65){\line(-4,-1){68}}
\put(118,65){\line(-4,-1){68}}
\put(158,65){\line(-4,-1){68}}
\put(178,65){\line(-4,-1){68}}
\put(102,65){\line(4,-1){68}}
\put(122,65){\line(4,-1){68}}
\put(162,65){\line(4,-1){68}}
\put(182,65){\line(4,-1){68}}

\put(300,23){\makebox(0,0){$(k^1_{A,B}+k^1_{C,D})\times(k^2_{A,B}+k^2_{C,D})$}}
\put(280,55){\makebox(0,0){$w_{(A\times B)+(C\times D)}$}}
\end{picture}

\begin{picture}(350,80)
\put(90,10){\makebox(0,0){$($}}
\put(100,10){\makebox(0,0){$A$}}
\put(110,10){\makebox(0,0){$+$}}
\put(120,10){\makebox(0,0){$C$}}
\put(130,10){\makebox(0,0){$)$}}
\put(140,10){\makebox(0,0){$\times$}}
\put(150,10){\makebox(0,0){$($}}
\put(160,10){\makebox(0,0){$B$}}
\put(170,10){\makebox(0,0){$+$}}
\put(180,10){\makebox(0,0){$D$}}
\put(190,10){\makebox(0,0){$)$}}

\put(10,40){\makebox(0,0){$($}}
\put(20,40){\makebox(0,0){$($}}
\put(30,40){\makebox(0,0){$A$}}
\put(40,40){\makebox(0,0){$+$}}
\put(50,40){\makebox(0,0){$C$}}
\put(60,40){\makebox(0,0){$)$}}
\put(70,40){\makebox(0,0){$\times$}}
\put(80,40){\makebox(0,0){$($}}
\put(90,40){\makebox(0,0){$B$}}
\put(100,40){\makebox(0,0){$+$}}
\put(110,40){\makebox(0,0){$D$}}
\put(120,40){\makebox(0,0){$)$}}
\put(130,40){\makebox(0,0){$)$}}
\put(140,40){\makebox(0,0){$+$}}
\put(150,40){\makebox(0,0){$($}}
\put(160,40){\makebox(0,0){$($}}
\put(170,40){\makebox(0,0){$A$}}
\put(180,40){\makebox(0,0){$+$}}
\put(190,40){\makebox(0,0){$C$}}
\put(200,40){\makebox(0,0){$)$}}
\put(210,40){\makebox(0,0){$\times$}}
\put(220,40){\makebox(0,0){$($}}
\put(230,40){\makebox(0,0){$B$}}
\put(240,40){\makebox(0,0){$+$}}
\put(250,40){\makebox(0,0){$D$}}
\put(260,40){\makebox(0,0){$)$}}
\put(270,40){\makebox(0,0){$)$}}

\put(90,70){\makebox(0,0){$($}}
\put(100,70){\makebox(0,0){$A$}}
\put(110,70){\makebox(0,0){$\times$}}
\put(120,70){\makebox(0,0){$B$}}
\put(130,70){\makebox(0,0){$)$}}
\put(140,70){\makebox(0,0){$+$}}
\put(150,70){\makebox(0,0){$($}}
\put(160,70){\makebox(0,0){$C$}}
\put(170,70){\makebox(0,0){$\times$}}
\put(180,70){\makebox(0,0){$D$}}
\put(190,70){\makebox(0,0){$)$}}

\put(100,65){\line(-4,-1){68}}
\put(120,65){\line(-3,-2){28}}
\put(160,65){\line(3,-2){28}}
\put(180,65){\line(4,-1){68}}

\put(30,35){\line(4,-1){68}}
\put(50,35){\line(4,-1){68}}
\put(90,35){\line(4,-1){68}}
\put(110,35){\line(4,-1){68}}
\put(170,35){\line(-4,-1){68}}
\put(190,35){\line(-4,-1){68}}
\put(230,35){\line(-4,-1){68}}
\put(250,35){\line(-4,-1){68}}

\put(280,23){\makebox(0,0){$m_{(A+C)\times (B+D)}$}}
\put(300,55){\makebox(0,0){$(l^1_{A,C}\times l^1_{B,D})+(l^2_{A,C}\times l^2_{B,D})$}}
\end{picture}
\end{center}

Each line in such a diagram stands for a family of parallel lines, one for
each propositional letter in the schemata $A$, $B$, $C$ and $D$.

In general, to verify whether for $f,g:A\rightarrow B$ in the language of $%
{\cal C}_{\times,+}$ we have $f=g$ in ${\cal C}_{\times,+}$ it is enough to
draw $G(f)$ and $G(g)$, and check whether they are equal, which is clearly a
finite task. So we have here an easy decision procedure for the equations of
${\cal C}_{\times,+}$.

\section{Maximality}

We shall now show that categories with binary products and sums, which we
call $\times\!,\!+$-categories, are maximal in the sense that if any
equation $f=g$ in the language of the free $\times\!,\!+$-category ${\cal C}%
_{\times,+}$ that doesn't hold in ${\cal C}_{\times,+}$ holds in a $%
\times\!,\!+$-category ${\cal B}$, then ${\cal B}$ is a preorder, i.e. a
category where all arrows with the same source and the same target are
equal. That $f=g$ holds in ${\cal B}$ means that it holds universally with
respect to objects. Namely, propositional letters in the indices of $f$ and $%
g$ are assumed to be variables, and $f=g$ holds iff it holds for every
assignment of objects to these variables. (This sort of universal holding is
quite natural in logic, and elsewhere in mathematics: it is usually taken
for granted when one says that a formula with variables ``holds''. In the
lambda calculus this universal holding is sometimes called ``typical
ambiguity''.) An analogous maximality is proved for ${\cal C}_{\times,%
\mbox{\scriptsize{\rm I}}}$ and ${\cal C}_{\times}$ (and, by duality, for $%
{\cal C}_{+,\mbox{\scriptsize{\rm O}}}$ and ${\cal C}_{+}$) in \cite{dos00},
and for cartesian closed categories in \cite{sim95} and \cite{dos00a}.

Suppose $A$ and $B$ are formulae of ${\cal P}_{\times,+}$ in which only $p$
occurs as a propositional letter. If for $f,g:A\rightarrow B$, we have $%
G(f)\neq G(g)$, then for some $x\in |A|$ and some $y\in |B|$ we have $%
(x,y)\in G(f)$ and $(x,y)\not\in  G(g)$, or vice versa. Suppose $(x,y)\in
G(f)$ and $(x,y)\not\in G(g)$. For every subformula $C$ of $A$ and every
formula $D$ let $A^C_D$ be the formula obtained from $A$ by replacing the
particular occurrence of the subformula $C$ by $D$. It is easy to see that
for every subformula $A_1+A_2$ of $A$ we have an arrow $h(l^i_{A_1,A_2})$ of
${\cal C}_{\times,+}$ built with $l^i_{A_1,A_2}$, identity arrows and the
operations on arrows $\times$ and $+$, such that $f\circ h(l^i_{A_1,A_2})$
and $g\circ h(l^i_{A_1,A_2})$ are of type $A^{A_1+A_2}_{A_i}\rightarrow B$.

We say that $x\in\omega$ {\em belongs to a subformula} $C$ of $A$ iff the $x$%
-th occurrence of propositional letters in $A$, counting from the left, is
in $C$. If $x$ happens to belong to $A_1$, we take care above to choose $%
h(l^1_{A_1,A_2})$, and if it belongs to $A_2$, we choose $h(l^2_{A_1,A_2})$.
If $x$ belongs to neither, we choose $h(l^i_{A_1,A_2})$ arbitrarily. By
repeated compositions of $f$ and $g$ with such $h(l^i_{A_1,A_2})$ arrows,
for every $+$ in $A$, we obtain two arrows $f^{\prime},g^{\prime}:p\times%
\ldots\times p\rightarrow B$ of ${\cal C}_{\times,+}$ such that parentheses
are somehow associated in $p\times\ldots\times p$, and for some $(z,y)\in
G(f^{\prime})$ we have $(z,y)\not\in G(g^{\prime})$. The formula $%
p\times\ldots\times p$ may be only $p$. We may further compose $f^{\prime}$
and $g^{\prime}$ with natural isomorphisms of the types $(C_1\times
C_2)\times C_3\rightarrow C_1\times(C_2\times C_3)$, $C_1\times (C_2\times
C_3)\rightarrow (C_1\times C_2)\times C_3$ and $C_1\times C_2\rightarrow
C_2\times C_1$, which are definable in ${\cal C}_{\times}$, and hence also
in ${\cal C}_{\times,+}$, in order to obtain two arrows $f^{\prime%
\prime},g^{\prime\prime}:p\times A^{\prime}\rightarrow B$ or $%
f^{\prime\prime},g^{\prime\prime}:p\rightarrow B$ such that $A^{\prime}$ is
of the form $p\times\ldots\times p$ with parentheses somehow associated, and
$(0,y)\in G(f^{\prime\prime})$ but $(0,y)\not\in G(g^{\prime\prime})$. The
functor $G$ maps the natural associativity and commutativity isomorphisms
into bijections.

By working dually on every $\times$ in $B$ using $h(k^i_{B_1,B_2})$, and by
composing perhaps further with natural associativity and commutativity
isomorphisms of $+$, we obtain two arrows $f^{\prime\prime\prime}$ and $%
g^{\prime\prime\prime}$ of ${\cal C}_{\times,+}$ of type $p\times
A^{\prime}\rightarrow p+B^{\prime}$ for $A^{\prime}$ of the form $%
p\times\ldots\times p$ and $B^{\prime}$ of the form $p+\ldots +p$, or of
type $p\times A^{\prime}\rightarrow p$, or of type $p\rightarrow
p+B^{\prime} $, such that $(0,0)\in G(f^{\prime\prime\prime})$ and $(0,0)
\not\in G(g^{\prime\prime\prime})$. (We cannot obtain that $%
f^{\prime\prime\prime}$ and $g^{\prime\prime\prime}$ are of type $%
p\rightarrow p$, since otherwise, by Cut Elimination, $g$ would not exist.)

With the help of $w_p$ we can define in ${\cal C}_{\times,+}$ the arrow $%
h^{\times}:p\rightarrow p\times\ldots\times p$ such that for every $x\in
|p\times\ldots\times p|$ we have $(0,x)\in G(h^{\times})$. We define
analogously with the help of $m_p$ the arrow $h^+:p+\ldots +p\rightarrow p$
such that for every $x\in |p+\ldots +p|$ we have $(x,0)\in G(h^{+})$.
(The arrows $h^{\times}$ and $h^+$ may be $\mbox{\bf 1}_p:p\rightarrow p$.)

If
$f^{\prime\prime\prime}$ and $g^{\prime\prime\prime}$ are of type $p\times
A^{\prime}\rightarrow p+B^{\prime}$,
let $f^{\dag},g^{\dag}:p\times p\rightarrow p+p$ be defined by
\[
\begin{array}{l}
f^{\dag}=_{def.}(\mbox{\bf 1}_p+h^+)\circ f^{\prime\prime\prime}\circ (%
\mbox{\bf 1}_p\times h^{\times}), \\[.05cm]
g^{\dag}=_{def.}(\mbox{\bf 1}_p+h^+)\circ g^{\prime\prime\prime}\circ (%
\mbox{\bf 1}_p\times h^{\times}).
\end{array}
\]
By Cut Elimination we have that $G(f^{\dag})$ and $G(g^{\dag})$
are singletons.
If $(1,0)$ or $(1,1)$ belong to $G(g^{\dag})$, then for
$f^*,g^*:p\times p\rightarrow p$ defined as $m_p\circ f^{\dag}$
and $m_p\circ g^{\dag}$, respectively,
we have $(0,0)\in G(f^*)$ and $(0,0)\not\in G(g^*)$.
If $(0,1)$ or $(1,1)$ belong to
$G(g^{\dag})$, then for
$f^*,g^*:p\rightarrow p+p$ defined as $f^{\dag}\circ w_p$
and $g^{\dag}\circ w_p$, respectively,
we have $(0,0)\in G(f^*)$ and $(0,0)\not\in G(g^*)$.

If $f^{\prime\prime\prime}$ and $g^{\prime\prime\prime}$ are of type $p\times
A^{\prime}\rightarrow p$, then for $f^*,g^*:p\times p\rightarrow p$ defined
as $f^{\prime\prime\prime}\circ (\mbox{\bf 1}_p\times h^{\times})$ and
$g^{\prime\prime\prime}\circ (\mbox{\bf 1}_p\times h^{\times})$, respectively,
we have $(0,0)\in G(f^*)$ and $(0,0)\not\in G(g^*)$.

If $f^{\prime\prime\prime}$ and $g^{\prime\prime\prime}$ are of type
$p\rightarrow p+B^{\prime}$, then for $f^*,g^*:p\rightarrow p+p$ defined as
$(\mbox{\bf 1}_p+h^+)\circ f^{\prime\prime\prime}$ and $(\mbox{\bf 1}%
_p+h^+)\circ g^{\prime\prime\prime}$, respectively, we have $(0,0)\in G(f^*)$
and $(0,0)\not\in G(g^*)$. In all that, we have
by Cut Elimination that $G(f^*)$ and $G(g^*)$ are singletons.

In cases where $f^*$ and $g^*$
are of type $p\times p\rightarrow p$, by
Cut Elimination, by the conditions on $G(f^*)$ and $G(g^*)$, and by the
functoriality of $G$, we obtain that $f^*=k^1_{p,p}$ and $g^*=k^2_{p,p}$. So
from $f=g$ we can derive $k^1_{p,p}=k^2_{p,p}$.
In cases where $f^*$ and $g^*$ are of type $p\rightarrow p+p$, by
Cut Elimination, by the conditions on
$G(f^*)$ and $G(g^*)$, and by the functoriality of $G$, we obtain that
$f^*=l^1_{p,p}$ and $g^*=l^2_{p,p}$. So from $f=g$ we can derive
$l^1_{p,p}=l^2_{p,p}$.
In any case, from $f=g$ we can derive $k^1_{p,p}=k^2_{p,p}$ or
$l^1_{p,p}=l^2_{p,p}$. If either of these two equations holds in a
$\times\!,\!+$-category ${\cal B}$, then ${\cal B}$ is a preorder. For
$h_1,h_2:C\rightarrow D$ in ${\cal B}$ we have
\[
\begin{array}{l}
k^1_{D,D}\circ\langle h_1,h_2\rangle = k^2_{D,D}\circ\langle h_1,h_2\rangle,
\;{\mbox{\rm or }} \\[.05cm]
[h_1,h_2]\circ l^1_{C,C}=[h_1,h_2]\circ l^2_{C,C},
\end{array}
\]
from which $h_1=h_2$ follows. (We have said that the holding of $%
k^1_{p,p}=k^2_{p,p}$ or $l^1_{p,p}=l^2_{p,p}$ is understood universally with
respect to objects of ${\cal B}$, so that we may replace $p$ by any object
of ${\cal B}$.)

It remains to remark that if for any $f,g:A\rightarrow B$ of ${\cal C}%
_{\times,+}$ we have that $f=g$ doesn't hold in ${\cal C}_{\times,+}$, then
by the Faithfulness of $G$ from ${\cal C}_{\times,+}$ we have $G(f)\neq G(g)$%
. If we take the substitution instances $f^{\prime}$ and $g^{\prime}$ of $f$
and $g$ obtained by replacing uniformly every propositional letter in $A$
and $B$ by $p$, then we obtain again that $G(f^{\prime})\neq G(g^{\prime})$.
If $f=g$ holds in a $\times\!,\!+$-category ${\cal B}$, then $%
f^{\prime}=g^{\prime}$ holds too, and hence ${\cal B}$ is a preorder, as we
have shown above.

This maximality result means that all equations in the language of ${\cal C}%
_{\times,+}$ that don't hold in ${\cal C}_{\times,+}$ can be derived from
each other with the help of the equations of ${\cal C}_{\times,+}$. This
result is analogous to the syntactic completeness of the classical
propositional calculus discovered by Bernays and Hilbert (see \cite{zac99},
p. 341), which is called Post Completeness. The Faithfulness of $G$ from $%
{\cal C}_{\times,+}$ amounts to a semantical completeness result (soundness
is provided by $G$ being a functor).

As a consequence of the maximality of $\times\!,\!+$-categories we obtain
that ${\cal C}_{\times,+}$ is a subcategory of ${\cal C}$. By Cut
Elimination we may conclude that it is a full subcategory of ${\cal C}$.
This means that the equations between terms of ${\cal C}$ make a
conservative extension of the equations of ${\cal C}_{\times,+}$; namely, if
an equation in the language of ${\cal C}_{\times,+}$ holds in ${\cal C}$,
then it holds already in ${\cal C}_{\times,+}$. This applies also to any
other category that is not a preorder whose equations extend those of ${\cal %
C}_{\times,+}$: for example, the free distributive bicartesian category
generated by ${\cal L}$, or the free bicartesian closed category generated
by ${\cal L}$.

So we may use the decision procedure provided by our Faithfulness of $G$
from ${\cal C}_{\times,+}$ in order to check equations of arrows in these
extensions of ${\cal C}_{\times,+}$, provided the terms of these arrows are
terms of ${\cal C}_{\times,+}$.

Bicartesian categories are not maximal in the sense in which $\times\!,\!+$%
-categories are. {\bf Set} is a bicartesian category that is not a preorder
in which $k^1_{\mbox{\scriptsize{\rm O}},\mbox{\scriptsize{\rm O}}}=k^2_{%
\mbox{\scriptsize{\rm O}},\mbox{\scriptsize{\rm O}}}$ holds, though this
equation does not hold in every bicartesian category. (Another example of
such an equation is $l_{\mbox{\scriptsize{\rm O}} \times A}\circ k^1_{%
\mbox{\scriptsize{\rm O}}, A}=\mbox{\bf 1}_{\mbox{\scriptsize{\rm O}} \times
A}$.) In a companion to this paper \cite{dos00c} we shall show that it is
enough to add $k^1_{\mbox{\scriptsize{\rm O}},\mbox{\scriptsize{\rm O}}%
}=k^2_{\mbox{\scriptsize{\rm O}},\mbox{\scriptsize{\rm O}}}$ to ${\cal C}%
_{\times,+,\mbox{\scriptsize{\rm O}}}$ to obtain coherence. We
shall show that we have only a restricted form of coherence when
we add to ${\cal C}$
this equation and $l^1_{\mbox{\scriptsize{\rm I}},\mbox{\scriptsize{\rm I}}}=l^2_{%
\mbox{\scriptsize{\rm I}},\mbox{\scriptsize{\rm I}}}$. However, this last
equation does not hold in {\bf Set}.

\end{document}